\documentclass{amsart}
 
 \usepackage{hyperref}    
  \usepackage{pdfsync}      
\usepackage{amsmath, amssymb, mathrsfs, amsfonts,  amsthm} 
\usepackage [all]{xy}

\newtheorem {theorem} {Theorem} [section]

\newtheorem {remark} [theorem] {Remark}
\newtheorem {rem/q}[theorem]{Remark/Question}
\numberwithin {equation} {section}

\author {Yasha Savelyev}
 \address{Centre de Recherches Math\'ematiques, Universit\'e de Montr\'eal, C.P.
6128, Succ. Centre-ville, Montr\'eal H3C 3J7, Qu\'ebec, Canada}
   \email{yasha.savelyev@gmail.com}
\title {Proof of the index conjecture in Hofer geometry}
\begin{document}
\begin{abstract} Let $\gamma$ be a non-degenerate Ustilovsky geodesic in $Ham
(M, \omega)$ generated by $H$. We give a simple proof  of a generalization of
the conjecture stated in \cite{virtmorse},  relating the Morse index of $
\gamma$, as a critical point of the Hofer length functional, with the Conley Zehnder index of
 the extremizers
of $H$, considered as  periodic orbits.  
\end{abstract}
\maketitle
\section {Introduction}
There has not been much study of the Morse index of geodesics for the Hofer
length functional on path spaces of the group of Hamiltonian diffeomorphisms $
\text {Ham}(M, \omega)$. Maybe this is because the problem of
Morse theory for the Hofer length functional seems completely hopeless. This is
possibly true to a large extent, however in \cite{virtmorse} we showed that doing
Morse theory for the Hofer length functional ``virtually'' can give some 
interesting results in symplectic topology. 
 
In the special case where $\gamma$
is an $S^1$-subgroup in $Ham (M, \omega) $ generated by a Morse Hamiltonian
$H$, a key point in \cite{virtmorse}
was using a
relationship of the Morse index of $ \gamma$ with the Conley-Zehnder index of
the linearized flow at the extremizers of $H$, in some special cases.
\begin{remark} We didn't use the words Conley-Zehnder index in \cite{virtmorse}, but
rather the index of a certain Cauchy-Riemann operator, but this could be directly related
to the above CZ index.
\end{remark}
 Indeed as a byproduct we arrived at the conjecture that the two indexes must
 coincide. A lower bound for the Morse index in terms of the Conley-Zehnder
 index was proved by Karshon-Slimowitz in \cite{Yael} by constructing a beautiful
  explicit local family of shortenings of $\gamma$.  Here we give a simple proof of the
 conjecture for more general Ustilovsky geodesics, using calculus of variations,
 already worked out in \cite{U} for the Hofer length functional.
   
In \cite{MTH} we use this coincidence to extend the virtual Morse theory
picture of \cite{virtmorse} from  special flag manifolds to general monotone
symplectic manifolds. 
 \subsection {Acknowledgements} I would like to thank
Leonid Polterovich who gave a crucial initial suggestion. Egor Shelukhin for
convincing me to consider the general case. And the anonymous referee for
carefully explaining an error regarding normalization in an earlier draft. This
paper was completed while the author was CRM-ISM postdoctoral fellow at CRM
Montreal.
\section {Statement and Proof}
\subsection {The group of Hamiltonian symplectomorphisms and Hofer metric}
\label{sec.hofer} Given a smooth function $H: M ^{2n} \times [0,1]
\to \mathbb{R}$,  there is  an associated time dependent
Hamiltonian vector field $X _{t}$, $0 \leq t \leq 1$,  defined by 
\begin{equation}  \label {equation.ham.flow} \omega (X _{t}, \cdot)=-dH _{t} (
\cdot).
\end{equation}
The vector field $X _{t}$ generates a path $\gamma: [0,1] \to \text
{Diff}(M)$, starting at $id$. Given such a path $\gamma$, its end point
$\gamma(1)$ is called a Hamiltonian symplectomorphism. The space of Hamiltonian symplectomorphisms
forms a group, denoted by $\text {Ham}(M, \omega)$.

In particular the path $\gamma$ above lies in $ \text {Ham}(M, \omega)$. It
is well-known 
that any smooth path $\gamma$ in $\text {Ham}(M, \omega)$ with
$\gamma (0)=id$ arises in this way (is generated by $H: M \times [0,1] \to
\mathbb{R}$ as above). Given  a general smooth path $\gamma$, the \emph{Hofer
length}, $L (\gamma)$ is defined by \begin{equation*}L (\gamma):= \int_0 ^{1} \max _{M} H _{t} ^{\gamma} 
-\min _{M} H ^{\gamma} _{t} dt,
\end{equation*}
where $H ^{\gamma}$ is a generating function for the path
$t \mapsto \gamma({0}) ^{-1} \gamma (t),$ $0\leq t \leq 1$.
The Hofer distance $\rho (\phi, \psi)$ is  defined by taking the
infinum of the Hofer length of paths from $\phi $ to $\psi$. We only mention
it, to emphasize that it is a deep and interesting theorem that the resulting
metric is non-degenerate, (cf. \cite{H, LM}). This gives $ \text {Ham}(M, \omega)$ the structure of a Finsler
manifold. 

We now consider $L$ as a functional on the space of paths in $ \text {Ham}(M,
\omega)$ starting at $id$ and ending at some fixed end point, denote this by
$\Omega \text {Ham}(M, \omega)$.
It is shown by Ustilovsky that $ \gamma$ is a smooth critical
point of $$L:  \Omega \text {Ham}(M, \omega) \to
\mathbb{R},$$ if there is a unique pair of points $x _{\max}$, $x _{\min} \in M$
maximizing, respectively minimizing the generating function $H ^{\gamma} _{t}$
at each moment $t$, and such that $H ^{\gamma} _{t}$ is Morse at $x _{\max}$, $x
_{\min}$. We shall call such a $\gamma$ \emph{Ustilovsky geodesic}. 

Consequently it makes sense to ask for the  Morse index of Ustilovsky geodesics,
(which might a priori be infinite.) Moreover, it is easy to see that 
$index _{L} (\gamma) = index _{L_+} (\gamma) + index _{L _-} (\gamma)$, 
where: \begin{equation} \label {eq.pos.length} L _{+} (\gamma):= \int_0 ^{1} \max
(H ^{\gamma} _{t}) dt, \end{equation} for $H ^{\gamma} _{t}$ in addition
normalized by the condition: 
\begin{equation} \label {eq.norm} \int _{M} H ^{\gamma} _{t} \cdot \omega
^{n}=0.
\end{equation}
The functional $L _{-}$ is  defined similarly as above. It will be the Morse
index of $\gamma$ with respect to $L _{+}$ that we compute.


Fix a small $\epsilon$, $0 <\epsilon <1$, s.t. the linearized flow (isotopy) at
$x _{\max}$ of $\gamma| _{ [0, \epsilon]}$ has no non trivial periodic orbits with positive
period.  Let us denote the periodic orbit of the isotopy
$\gamma| _{ [0, \epsilon]}$ associated to $x _{\max}$ by $x _{\max,0}$, and likewise the periodic orbit of
the isotopy $\gamma| _{ [0,1]}$ associated to $x _{\max}$ by $x _{\max,1}$. We
will say that $\gamma$ is \emph{non-degenerate} if $x _{\max, 1}$ is
non-degenerate in the sense of Floer theory, in other words the time 1 linearized flow at $x
_{\max}$ has no non-trivial time 1 perioidic orbits.

 \begin{theorem} For $\gamma$ a non-degenerate Ustilovsky geodesics as above,
 the Morse index of $ \gamma$ with respect to $L _{+}$ is
\begin{equation} \label {eq.difference} 
|CZ (x _{\max,1}) - CZ (x _{\max,0})|.
\end{equation}
\end{theorem}
\begin{remark}  The above expression is independent of any choices of normalization of $CZ$
index appearing in literature. Moreover it is precisely the
index of the real linear Cauchy Riemann operator on which the conjecture is
based in \cite{virtmorse}. A better way to understand this coincidence is
outlined in Section 1.3 of that paper.
\end{remark}

\begin{proof} 
The Morse index theorem \cite{Morse.index} cannot be directly applied to $$L _{+} = \int
_{0} ^{1} L (\dot{\gamma} (t), \gamma (t)) dt,$$ $L (\dot {\gamma} (t), \gamma
(t)) = \max _{M} H _{t}$, for $H _{t} = \dot {\gamma} (t) \in T _{\gamma (t)}
\text {Ham}(M, \omega) \equiv C ^{\infty} _{norm} (M)$, with the latter being
smooth functions normalized to have zero mean, \eqref{eq.norm}.
This is because it clearly does not satisfy the  Legendre condition that $
\frac{d^2}{d \xi ^{2}} L (\dot {\gamma} (t), \gamma
(t)) >0$, for every variation $\xi$ of $\dot { \gamma}$, (for every $t$). 
However Ustilovsky shows that there is a related functional (actually a
quadratic form) $\mathcal {L}_+$ on the vector space $\Omega_0 T
_{x_{\max}} M$ (based loop space at 0 on the tangent space). 
With the Hessian of $L_+$ at $\gamma$ coinciding with the
Hessian of $ \mathcal {L} _{+}$ at $0$, and to which the Morse theorem does
apply. This is beautiful, but we refer the reader to \cite{U} and
\cite[Section 12.4]{polterov.geom} for further details.

The Morse theorem gives us the following procedure for the
calculation of the Morse index of $\gamma$ with respect to $ \mathcal {L}_+$.
Denote by $\gamma _{\tau}$ the restriction of $\gamma$ to $ [0, \tau] \subset [0,1]$. Then
$index (\gamma _{\tau})$ is a locally constant, lower semi-continuous function
in $ \tau$, and jumps at a discrete set of $ \tau _{i} \in (0,1)$ called
conjugate times. The value of the jump  
$mult ( \tau_i)$ is the dimension of the solution space of the associated Jacobi
equation. Informally speaking this is dimension of the space of infinitesimal
variations of $\gamma _{\tau}$ through extremals with the same endpoints. And a
point $\tau \in (0, 1]$ is defined to be a conjugate time if this dimension is
non zero.
 
 In the case of the functional $ \mathcal {L}_+$, it is shown in \cite{U} that
 $\tau_0 \in (0,1]$ is a conjugate time if and only if the time $\tau_0$
 linearized flow of $H$ at the extremizer $x _{\max}$ of $H$ has periodic orbits, and the
 multiplicity  $mult ( \tau _{0})$ is the dimension of the space of these
 periodic orbits. 
 


To keep notation simple, let us denote by $\gamma _{\max}$ the
restriction of the linearization of $\gamma$ at $x _{\max}$ to $ [\epsilon, 1]$. 
We will use the
construction of Maslov and Conley-Zehnder index given in
\cite{Maslov.for.paths}. For the normalizations used in \cite{Maslov.for.paths} we show that the absolute value 
of the Conley-Zehnder index for the path $\gamma _{\max}$ is exactly the Morse index of $\gamma$ for $L _{+}$, from which the
statement of the theorem immediately follows by additivity of the
Conley-Zehnder/Maslov index with respect to concatenations (and with respect to those normalizations). 

Note first that $\gamma _{\max}$ has a crossing at
$\tau_0 \in (0,1)$ with the Maslov cycle  if and only if $\gamma _{\max}
(\tau_0)$ has 1-eigenvectors, i.e. if and only if  $\tau_0$ is a
conjugate time. Moreover, the dimension of the intersection $I _{\tau_0}$ of the
diagonal $\Lambda \subset \mathbb{R} ^{2n} \times \mathbb{R} ^{2n}$ with the graph $Gr
(\gamma _{\max} (\tau_0)) = \{(z, \gamma _{\max} (\tau_0) z)| z \in \mathbb{R}
^{2n}\}$, is exactly the multiplicity of $\tau _{0}$. The crossing form $Q$ at
$\tau _{0}$ can then be identified with the Hessian of $H ^{\gamma} _{\tau
_{0}}$ at $x _{\max}$, which follows by \cite[Remark 5.4]{Maslov.for.paths}.
Since this is non-degenerate by assumption, all the crossings are regular. 
 Our conventions are
\begin{align*} \omega (X _{H}, \cdot) = -dH (\cdot)\\  
\omega (\cdot, J \cdot)>0.
\end{align*} 
Consequently the crossing form is negative definite, and so is negative definite
on $I _{\tau_0}$. So the signature of $Q$ on $I _{\tau_0}$ (number of positive
minus number of negative eigenvalues) is just the $- mult (\tau_0)$. The
Conley-Zehnder index of $\gamma _{\max}$ is then the sum over conjugate times
$\tau _{i}$ of $-mult (\tau _{i})$. Consequently Morse index of $\gamma$, is $|CZ (\gamma _{\max})|$.

\end{proof}   
\bibliographystyle{siam}

\end{document}